\documentclass[letterpaper, 10 pt, conference]{ieeeconf}  % Comment this line out if you need a4paper
\usepackage{xcolor}
\usepackage{graphicx}
\usepackage{mathtools}
\usepackage{tabularx} 
\usepackage{booktabs}
\usepackage{amsmath,amssymb,amsfonts}
\usepackage{graphicx}
\usepackage{amssymb}
\usepackage{eqnarray} 
\usepackage{amsmath}
\usepackage{mathalfa} 
\usepackage{amsmath}
\usepackage{stfloats}
\usepackage{amssymb}
\usepackage{dsfont}
\usepackage{xcolor}
\usepackage{bbm}
\usepackage{float}
\usepackage{MnSymbol}
\usepackage{mathrsfs}
\usepackage{lscape}
\usepackage{longtable}
\usepackage{rotating}
\usepackage{multirow}
\usepackage{color}
\usepackage{xcolor}
\usepackage{url}
\usepackage{subfigure}
\usepackage{rotating}
\DeclareMathOperator{\support}{supp}

\newcommand{\sB}{\mathcal{B}}
\newcommand{\sC}{\mathcal{C}}
\newcommand{\sD}{\mathcal{D}}

\newcommand{\sR}{\mathcal{R}}

\newcommand{\sX}{\mathcal{X}}
\newcommand{\sU}{\mathcal{U}}

\newcommand{\R}{\mathbb{R}}
\newcommand{\E}{\mathbb{E}}
\newcommand{\Z}{\mathbb{Z}}
\newcommand{\K}{\mathbb{K}}

\newcommand{\N}{\mathbb{N}}

\newcommand{\bv}[1]{\mathbf{#1}}
\newcommand{\statecost}[1]{c_{#1}^{(x)}\left(\bv{x}_{#1}\right)}
\newcommand{\actioncost}[1]{c_{#1}^{(u)}\left(\bv{u}_{#1}\right)}
\newcommand{\statecostexpectation}[1]{c_{#1}^{(x)}\left(\bv{X}_{#1}\right)}

\newcommand{\Vx}[2]{V_{\alpha}\left(\bv{x}_{k-1},\bv{u}_k\right)}
\newcommand{\Wx}[2]{W_{\alpha}\left(\bv{x}_{k-1},\bv{u}_k\right)}
\newcommand{\Vxtilde}[2]{\tilde{V}_{\alpha}\left(\bv{x}_{k-1},\bv{u}_k\right)}
\newcommand{\Mx}[2]{M\left(\bv{x}_{k-1},\bv{u}_k\right)}

\newcommand{\ratio}[2]{r_{{#1}\mid{#2}}}

\newcommand{\optimalratio}[2]{r^{\ast}_{{#1}\mid{#2}}}

\newcommand{\jointxu}[2]{p_{{#1}}\left(\bv{x}_{{#1}},\bv{u}_{{#1}}\mid \bv{x}_{{#2}} \right)}
\newcommand{\shortjointxu}[2]{p_{{#1}\mid{#2}}}
\newcommand{\refjointxu}[2]{q_{{#1}}\left(\bv{x}_{{#1}},\bv{u}_{{#1}}\mid \bv{x}_{{#2}} \right)}
\newcommand{\shortrefjointxu}[2]{q_{{#1}\mid{#2}}}
\newcommand{\shorttwistedrefjointxu}[2]{\tilde{q}_{{#1}\mid{#2}}}
\newcommand{\twistedrefjointxu}[2]{\tilde{q}_{{#1}}\left(\bv{x}_{{#1}},\bv{u}_{{#1}}\mid \bv{x}_{{#2}} \right)}

\newcommand{\plant}[2]{p_{{#1}} \left(\bv{x}_{{#1}}\mid \bv{x}_{{#2}}, \bv{u}_{{#1}} \right)}
\newcommand{\shortplant}[2]{p_{{#1}\mid{#2}}^{(x)}}
\newcommand{\refplant}[2]{q_{{#1}} \left(\bv{x}_{{#1}}\mid \bv{x}_{{#2}}, \bv{u}_{{#1}} \right)}
\newcommand{\shortrefplant}[2]{q_{{#1}\mid{#2}}^{(x)}}
\newcommand{\nominalplant}[2]{\bar{q}_{{#1}} \left(\bv{x}_{{#1}}\mid \bv{x}_{{#2}}, \bv{u}_{{#1}} \right)}
\newcommand{\shortnominalplant}[2]{\bar{q}_{{#1}\mid{#2}}^{(x)}}

\newcommand{\policy}[2]{p_{{#1}}\left(\bv{u}_{{#1}}\mid \bv{x}_{{#2}} \right)}
\newcommand{\shortpolicy}[2]{p_{{#1}\mid{#2}}^{(u)}}
\newcommand{\refpolicy}[2]{q_{{#1}}\left(\bv{u}_{{#1}}\mid \bv{x}_{{#2}} \right)}
\newcommand{\shortrefpolicy}[2]{q_{{#1}\mid{#2}}^{(u)}}

\newcommand{\optimalpolicy}[2]{p^{\star}_{{#1}} \left(\bv{u}_{{#1}}\mid \bv{x}_{{#2}} \right)}
\newcommand{\shortoptimalpolicy}[2]{p_{{#1}\mid{#2}}^{(u),\star}}

\newcommand{\ball}[2]{\sB_{\eta}\left(\shortnominalplant{#1}{#2}\right)}

\newcommand{\KL}{\text{KL}}

\newcommand{\radius}[2]{\eta_{#1}\left(\bv{x}_{#2},\bv{u}_{#1}\right)}

\newcommand{\DKL}[2]{{D}_{\KL}\left(#1\mid \mid #2 \right)}

\newtheorem{theorem}{Theorem}[section]

\newtheorem{definition}{Definition}[section]

\newtheorem{lemma}{Lemma}[section]

\newtheorem{problem}{Problem}[section]

\newtheorem{remark}{Remark}[section]

\newtheorem{assumption}{Assumption}[section]

\usepackage{dsfont}

\DeclareMathOperator{\argmin}{\operatorname{arg min}}

\usepackage{bbm}  % in the preamble

\usepackage{algorithm}
\usepackage{algorithmic}
\newenvironment{proofsketch}[1][Sketch of the Proof]{\textit{#1:} }{\ \rule{0.5em}{0.5em} \vspace{1ex}}

\IEEEoverridecommandlockouts                              % This command is only needed if 
\title{\LARGE \bf
Free-Energy Minimizing Policies Under Generative Model Ambiguity}

\author{Arash Shafiei$^{1,\ast}$, Caio César Graciani Rodrigues$^{2,\ast}$ and Giovanni Russo$^{3}$ % <-this % stops a space
\thanks{$^{1}$ A. Shafiei {\tt\small (shafie.allahkaram@gmail.com)} is with Czech Technical University  $^{2}$ C. C. Graciani Rodrigues {\tt\small (c.graciani@ssmeridonale.it)} is with Scuola Superiore Meridionale,  Italy, Modelling and Engineering Risk and Complexity Department, $^{3}$ G. Russo {\tt\small (giovarusso@unisa.it)}  is with University of Salerno, Italy, Department of Information and Electrical Engineering and Applied Mathematics. $^{\ast}$ Joint first authors %
}}

\begin{document}

\maketitle
\thispagestyle{empty}
\pagestyle{empty}

%%%%%%%%%%%%%%%%%%%%%%%%%%%%%%%%%%%%%%%%%%%%%%%%%%%%%%%%%%%%%%%%%%%%%%%%%%%%%%%%
\begin{abstract}
We present a variational free-energy formulation for distributionally robust decision-making with  ambiguity in the generative model. The formulation, related to a broad range of learning and control frameworks, yields a minimax optimal control problem where maximization is  over an uncertainty set that represents ambiguities. We prove that computing the optimal policy requires  solving a non-convex minimization problem and propose an algorithm with convergence guarantees to find the solution. The effectiveness of our results is illustrated via simulations on a pendulum swing-up problem.
\end{abstract}

%%%%%%%%%%%%%%%%%%%%%%%%%%%%%%%%%%%%%%%%%%%%%%%%%%%%%%%%%%%%%%%%%%%%%%%%%%%%%%%%
\section{Introduction}

The minimization of the variational free energy is central to a broad range of decision-making frameworks across learning, control and neuroscience~\cite{PG_MR_RW:14,FRISTON2009293}. These functionals consist of an expected cost and a term that quantifies the discrepancy from a generative model. Despite the success of these frameworks, free energy minimizing  policies can be brittle to ambiguities in the generative model, arising from, e.g., mis-specifications. In this context, a key challenge is guaranteeing robustness of free energy minimizing policies against possibly adversarial uncertainties in the generative model. Motivated by this, we present a distributionally robust variational free-energy formulation for policy computation under generative model ambiguity. The formulation yields an optimal control problem involving a non-convex distributionally robust optimization (DRO) problem and we introduce an algorithm to compute its optimal solution.\\

\noindent{\bf Related Works And Contributions.} Since Scarf's pioneering work, DRO has emerged as a central topic in optimization, gaining attention across machine learning and control, see, e.g.,~\cite{pmlr-v108-kirschner20a,BVP_DK_PG_MM:15,VanParys2021}.
%
%~\cite{pmlr-v258-dirren25a,pantazis2025datadrivenwassersteindistributionallyrobust} eu adicionei
%pmlr-v108-kirschner20a,BVP_DK_PG_MM:15
%Rahimian2022%VanParys2021
%buerger2012distributed
As also surveyed in~\cite{Lin2022,Xu2012,pmlr-v258-dirren25a}, these problems often feature a cost that is linear in at least one of the decision variables with extensions to non-linear, convex, settings rapidly emerging~\cite{Sheriff2024}. Popular frameworks such as distributionally robust versions of imitation learning based on maximum entropy~\cite{NEURIPS2021_cc8090c4}, Q-learning~\cite{pmlr-v162-liu22a} and path integral control~\cite{10644179,7862201} tackle policy computation problems framed via DRO where the cost is linear in the decision variables for the inner optimization problem.
In the context of control, contributions include, e.g., the development of distributionally robust linear quadratic controllers~\cite{10.5555/3666122.3666939} and model predictive control algorithms~\cite{10753052} with results benefiting from this linearity assumption. The minimization of free energy functionals arises across, e.g., MaxEnt, MaxDiff learning~\cite{eysenbach2022maximum,Berrueta2024}, KL Control~\cite{ET:09,PG_MR_RW:14}, Shr\"odinger bridges~\cite{Chen2015,doi:10.1137/20M1320195}, optimal~\cite{Srivastava2023} and inverse optimal control schemes based on maximum likelihood~\cite{KD-ET:10,EG_HJ_CDV_GR:25}. Related functionals also appear in the context of fully probabilistic design~\cite{MK:96,DG-GR:22}. We are, however, not aware of frameworks that tackle the problem of robustly minimizing this functional against possibly adversarial uncertainties in the generative model. In~\cite{shafiei2025robustdecisionmakingfreeenergy}, the problem of computing distributionally robust free energy minimizing policies has been considered and the resulting distributionally robust free energy principle has been integrated with learning frameworks in~\cite{jesawada2026learningbasedrobustcontrolunifying}. However, in~\cite{shafiei2025robustdecisionmakingfreeenergy,jesawada2026learningbasedrobustcontrolunifying} robustness is guaranteed against ambiguities in the environment. In this work, it is the generative model to be ambiguous and this, as we shall see, leads to a non-convex problem that cannot be studied with the tools in~\cite{shafiei2025robustdecisionmakingfreeenergy,jesawada2026learningbasedrobustcontrolunifying}.

Our key technical contributions are as follows. We:
\begin{enumerate}
    \item introduce a distributionally robust free energy formulation for policy computation. This yields a bi-level optimization approach to solve the  problem. The inner optimization step involves solving a non-convex, constrained and infinite-dimensional program. We obtain a difference-of-convex (DC) functions reformulation of the problem and characterize its main properties;
   % In contrast, ~\cite{kent2021frank,kent2021modified} formulate their inner problems via local linear minimization or Wasserstein-ball duality, and~\cite{Sheriff2024} employs a G-derivative framework without any DC or entropy-based structure;
   \item adapt a Frank-Wolfe algorithm from~\cite{millan2023frank} to address the optimization. The algorithm, which was originally introduced for vector spaces, is extended to infinite-dimensional spaces to tackle the settings required by our formulation.
   Although the Frank-Wolfe algorithms in~\cite{kent2021frank,kent2021modified} also address infinite-dimensional problems, they apply  to unconstrained settings, whereas our formulation is constrained;
    \item characterize convergence of the algorithm. Although convergence of the Frank-Wolfe algorithm has been  characterized for unconstrained problems in the probability space~\cite{kent2021frank,kent2021modified}, to tackle our formulation we need to guarantee algorithm convergence in the space of measurable functions with nonlinear constraints. To the best of our knowledge, this broad setting has not been tackled in the literature;
\item illustrate the effectiveness of our results via a numerical example.
%\item\color{red}{A.SH: Although the problems studied in~\cite{kent2021frank,kent2021modified} are formulated in infinite-dimensional settings, they remain within the class of unconstrained optimization problems. Similarly, in~\cite{Sheriff2024}, the authors consider a finite-dimensional min–max problem where the decision variable lies in a finite space and the constraint set corresponds to the entire space of probability measures. Consequently, the Frank–Wolfe (FW) methods proposed in those works are not directly applicable to our setting, which involves constrained optimization in infinite-dimensional functional spaces. In contrast, in our formulation the inner maximization (robustness) step leads to a non-convex optimization problem that cannot be addressed using conventional DRO or convex duality techniques. The main algorithmic novelty of this paper is therefore the development of a Frank–Wolfe–based method tailored for non-convex, infinite-dimensional, difference-of-convex (DC) constrained optimization problems arising from the proposed distributionally robust free-energy framework. This algorithm extends the classical Frank–Wolfe scheme to spaces of measurable functions, enabling optimization over probability densities and likelihood ratios, which is essential for explicitly accounting for ambiguity in generative models.  }
\end{enumerate}

%\textcolor{blue}{**** Caio, I commented your text to get an idea of the length of the paper. Please revise contributions accordingly. I went back to Arash notes but I may be wrong in differentiating from~\cite{kent2021frank,kent2021modified}.

%Embed all the novelties directly in the statement of contributions following the style I am using.

\section{Mathematical background}

Sets and operators are in {\em calligraphic} characters and vectors in {\bf bold}.
The symbol $\coloneqq$ denotes a definition and we let $\mathbb{K}$ be either $\R$ or $\Z$.
A random variable is denoted by $\bv{W}$ and its realization is $\bv{w}$. We denote by $p(\mathbf{w})$ the \textit{probability mass function} (pmf) when $\mathbf{W}$ is a discrete random variable, and the \textit{probability density function} (pdf) when $\mathbf{W}$ is continuous.
The convex subset of pdfs (pmfs) is denoted by $\sD$.
The expectation of a function $\mathbf{h}(\cdot)$ with respect to a continuous random variable $\mathbf{W}$ with density $p(\mathbf{w})$ is denoted by $\mathbb{E}_{p}[\mathbf{h}(\mathbf{w})] \coloneqq \int_{\operatorname{supp} p} \mathbf{h}(\mathbf{w})\, p(\mathbf{w})\, d\mathbf{w}$, where $\operatorname{supp} p$ denotes the (compact) support of $p(\mathbf{w})$. In the case where $\mathbf{W}$ is a discrete random variable, the integral is replaced by a summation, yielding $\mathbb{E}_{p}[\mathbf{h}(\mathbf{w})] \coloneqq \sum_{\mathbf{w}} \mathbf{h}(\mathbf{w})\, p(\mathbf{w})$. Whenever the context is clear, we omit the subscript in the expectation.
The joint pmf of $\mathbf{w}_1$ and $\mathbf{w}_2$ is denoted by $p(\mathbf{w}_1,\mathbf{w}_2)$ and the conditional pmf of $\mathbf{w}_1$ with respect to (w.r.t.) $\mathbf{w}_2$ is $p\left( \mathbf{w}_1\mid  \mathbf{w}_2 \right)$.  Given two pdfs, $p(\bv{w})$ and $q(\bv{w})$, we say that $p(\bv{w})$ is absolutely continuous with respect to (w.r.t.) $q(\bv{w})$ if $\support p \subseteq \support q$.  We denote this by writing $p \ll q$. %Countable sets are denoted by $\lbrace w_k \rbrace_{k_1:k_n}$, where $w_k$ is the generic set element, $k_1$ ($k_n$) is the index of the first (last) element and  $k_1:k_n$ is the set of consecutive integers between (including) $k_1$ and $k_n$.
%\textcolor{blue}{Let $(\Omega,\mathcal{F},\mu)$ denote the measurable space with $\Omega \subseteq \mathbb{R}^n$, $\mathcal{F} = \mathcal{B}(\Omega)$ the Borel $\sigma$-algebra and the Lebesgue measure.}
Also, let $\Omega\subseteq\K^n$, then $L^1(\Omega)$ stands for the space of measurable functions $f \colon \Omega \to \mathbb{R}$, such that: (i) 
$\|f\|_{L^1} = \int_\Omega |f| \, d\mu < \infty$, with Lebesgue measure $\mu$ if $\K^n$ is $\R^n$; (ii) $\|f\|_{L^1} = \sum_{x \in \Omega} |f(x)| < \infty$ if $\K^n$ is $\Z^n$.
We use the shorthand notation {\em a.s.} for {\em almost surely}.
The Kullback--Leibler divergence between distributions $p(\mathbf{w})$ and $q(\mathbf{w})$ is $    D_{\mathrm{KL}}(p \| q) \coloneqq \int_{\operatorname{supp} p} p(\mathbf{w}) \log \left( \frac{p(\mathbf{w})}{q(\mathbf{w})} \right) d\mathbf{w}$, if $\mathbf{W}$ is continuous. The integral
is replaced by a summation when $p(\mathbf{w})$ is a pmf.
The divergence -- measuring the discrepancy of the pair $p(\mathbf{w})$, $q(\mathbf{w})$ -- is non-negative and $D_{\mathrm{KL}}(p \| q) = 0$ if and only if $p = q$ almost everywhere \cite{SK_RL:51}. Finally, we recall the chain rule for the KL divergence~\cite{TC_JT:06}
\begin{lemma}\label{lem:chain_rule}
let $\bv{W}_1$ and $\bv{W}_2$ random variables, $f(\bv{w}_1,\bv{w}_2)$, $g(\bv{w}_1,\bv{w}_2)$  joint pdfs/pmfs and $f(\bv{w}_2\mid\bv{w}_1)$, $g(\bv{w}_2\mid\bv{w}_1)$ conditional pdfs/pmfs. Then: $\DKL{f(\bv{w}_1,\bv{w}_2)}{g(\bv{w}_1,\bv{w}_2)} = \DKL{f(\bv{w}_1)}{g(\bv{w}_1)} + \E_{f(\bv{w}_1)}\left[\DKL{f(\bv{w}_2\mid\bv{w}_1)}{g(\bv{w}_2\mid\bv{w}_1)}\right]$.
\end{lemma}

\subsection{Variational Free Energy}
The variational free energy, see, e.g.,~\cite{9363495}, is defined as
\begin{equation}\label{eqn:variational_free_energy}
\mathcal{F}(p(\mathbf{w}),q(\mathbf{w})) := D_{\mathrm{KL}}\left(p(\mathbf{w}) \mid \mid q(\mathbf{w})\right) + \E_p\left[l(\mathbf{W})\right],
\end{equation}
where $p(\mathbf{w})$ and $q(\mathbf{w})$ -- termed as generative (time-series) model -- are pdfs/pmfs and $l(\cdot)$ is a loss function. The following result can be found under a wide range of different technical conditions in, e.g.,~\cite{PG_MR_RW:14,Cammardella2023,ET:09}.
%ET:09
%DG_GR:22
\begin{lemma}\label{lem:free_energy_optimization}
given $q(\mathbf{w})$, consider 
$$
\underset{p(\mathbf{w})\in\sD}{\min} \mathcal{F}(p(\mathbf{w}),q(\mathbf{w})),
$$
with $\mathcal{F}(\cdot,\cdot)$ defined in~\eqref{eqn:variational_free_energy}. Then:
\begin{itemize}
    \item the problem is convex;
    \item the optimal solution is $p^{\ast}(\mathbf{w}) =  \frac{q(\mathbf{w})\exp\left(-l(\bv{w})\right)}{\int q(\mathbf{w})\exp\left(-l(\bv{w})\right)d\bv{w}}$;
    \item the  optimal value is $-\log\int q(\mathbf{w})\exp\left(-l(\bv{w})\right)d\bv{w}$.
\end{itemize}
\end{lemma}
\begin{remark}
for discrete variables, the integral in the optimal solution of Lemma~\ref{lem:free_energy_optimization} is replaced with a sum. The optimal value  arises in the context of, e.g., statistical physics~\cite[Chapter 20]{Coolen2005}.
\end{remark}
%\begin{remark}
%The KL divergence in~\eqref{eqn:variational_free_energy} can be thought of as a regularizer -- biasing the optimal solution of the problem in Lemma~\ref{lem:free_energy_optimization} towards $q(\bv{w})$ -- and a regularization term can be included in the formulation. 
%\end{remark}

\section{Problem statement}

Let {$\bv{X}_k\in\sX\subseteq\mathbb{K}^n$ be the state and $\bv{U}_k\in\sU\subseteq\mathbb{K}^m$} the control input at time-step $k$.  The possibly nonlinear and stochastic dynamics for the system is $\plant{k}{k-1}$ and we use the shorthand notation $\shortplant{k}{k-1}$ to denote this probability. We also let $\shortjointxu{k}{k-1} := \jointxu{k}{k-1} = \plant{k}{k-1}\policy{k}{k-1}$, where $\policy{k}{k-1}$ is a randomized policy. We use $\shortpolicy{k}{k-1}$ to denote this policy.

\begin{remark}
$\plant{k}{k-1}$ describes the behavior of the system at time-step $k$ given $\bv{x}_{k-1}$ and  $\bv{u}_k\sim\policy{k}{k-1}$. The joint $\jointxu{k}{k-1}$ captures in probabilistic terms the closed loop interactions. With the control problem, formalized next, we aim to design the agent policy.
\end{remark}

The policy is computed by minimizing, in the space of policies, the worst case variational free energy across a set capturing ambiguity on  a given $\shortnominalplant{k}{k-1}:=\nominalplant{k}{k-1}\in\sD$ serving as nominal generative model in our formulation. To formally state the problem, we first give the following:
\begin{definition}\label{def:ambiguity_set} 
let $\shortrefplant{k}{k-1}\coloneqq\refplant{k}{k-1}\in\sD$ and $\radius{k}{k-1}$ be a positive and bounded function. The ambiguity set associated with $\shortnominalplant{k}{k-1} \in \sD$ is
\begin{multline}\label{eqn:ambiguity_set}
\ball{k}{k-1} \coloneqq \left\{ \shortrefplant{k}{k-1} \in \sD : 
\DKL{\shortrefplant{k}{k-1}}{\shortnominalplant{k}{k-1}} \right. \\ \le \radius{k}{k-1}, 
\left. \support{\shortplant{k}{k-1}} \subseteq \support{\shortrefplant{k}{k-1}} \right\}.
\end{multline}
\end{definition}
\begin{remark}
$\ball{k}{k-1}$ is the set of $\refplant{k}{k-1}\in\sD$ that: (i) have statistical complexity, i.e., KL-divergence, of at most $\radius{k}{k-1}$ from $\nominalplant{k}{k-1}$; (ii) have a support that covers the support of $\plant{k}{k-1}$.
\end{remark}

We are now ready to state our control problem:
\begin{problem}\label{prob:main}
let:
\begin{enumerate}
\item $\shortplant{k}{k-1}\coloneqq\plant{k}{k-1}$;
\item {$\shortrefjointxu{k}{k-1}\coloneqq\refjointxu{k}{k-1} = \shortrefplant{k}{k-1}\shortrefpolicy{k}{k-1}$, with $\shortrefplant{k}{k-1} \coloneqq \refplant{k}{k-1}$ and $\shortrefpolicy{k}{k-1} \coloneqq \refpolicy{k}{k-1}$;}
\end{enumerate}
find $\shortoptimalpolicy{k}{k-1}\coloneqq\optimalpolicy{k}{k-1}$ such that
\begin{equation}\label{eqn:main_problem}
    \shortoptimalpolicy{k}{k-1} \in \text{arg } \underset{\shortpolicy{k}{k-1}\in\sD}{\min} \ \  \underset{\shortrefplant{k}{k-1}\in\ball{k}{k-1}}{\max} \ \mathcal{F} \left(\shortjointxu{k}{k-1},\shortrefjointxu{k}{k-1}\right),
\end{equation}
where:
%$\jointxu{k}{k-1}\coloneqq\shortplant{k}{k-1}\shortpolicy{k}{k-1}$ and
(i) $\mathcal{F}(\cdot,\cdot)$ is defined as in~\eqref{eqn:variational_free_energy} with  $l(\bv{x}_{k},\bv{u}_k) = \statecost{k} + \actioncost{k}$; (ii) $\ball{k}{k-1}$ is in Definition~\ref{def:ambiguity_set}.
\end{problem}
\begin{remark}
in Problem~\ref{prob:main}, the loss function in the free energy is a cost, which is possibly non-convex in states/actions. Minimizing the KL divergence term amounts at minimizing the discrepancy between $\jointxu{k}{k-1}$ and  the generative model. This model does not need to be a time-series model and can simply capture a desired configuration. In this case, by minimizing the KL divergence, the policy  minimizes the error from the desired configuration.
\end{remark}
\begin{remark}
in Problem~\ref{prob:main}, the constraint of the inner maximization problem captures the ambiguity around the generative model $\refplant{k}{k-1}$, which is not available.
\end{remark}
Following~\cite{shafiei2025robustdecisionmakingfreeenergy}, before introducing the results we briefly relate our Problem~\ref{prob:main} with other decision-making frameworks.

\subsection{Problem~\ref{prob:main}
And Other Decision-Making Frameworks}

The objective of the problem in~\eqref{eqn:main_problem} can be written as
\begin{equation}\label{eqn:KL_divergence_min}
\DKL{\shortjointxu{k}{k-1}}{\shorttwistedrefjointxu{k}{k-1}}-\log \left( Z(\shortrefjointxu{k}{k-1}) \right) , 
\end{equation}
where
\begin{align}\nonumber
& \shorttwistedrefjointxu{k}{k-1} \coloneqq\twistedrefjointxu{k}{k-1} \\
& = \frac{\refjointxu{k}{k-1} \exp\left(-\statecost{k}-\actioncost{k}\right)}{Z(\shortrefjointxu{k}{k-1})}
\end{align}
with
\begin{align}\nonumber
&Z(\shortrefjointxu{k}{k-1}) \\
&\coloneqq \int_{\sU} \int_{\sX} \shortrefjointxu{k}{k-1} \exp\left(-\statecost{k}-\actioncost{k}\right) d \bv{x}_k d \bv{u}_k    
\end{align}
being the normalizer. Differently from~\cite{shafiei2025robustdecisionmakingfreeenergy} the normalizer depends on the decision variables of the optimization problem and therefore it cannot be dropped from the optimization. However, when there is no ambiguity, the ambiguity constraint is relaxed and the second term in~\eqref{eqn:KL_divergence_min} can be dropped from the resulting minimization problem in the policy space. 

Next, we highlight some of the connections drawn in~\cite{shafiei2025robustdecisionmakingfreeenergy} between the relaxed minimization problem in the policy space and other frameworks; see~\cite{shafiei2025robustdecisionmakingfreeenergy} for details. In the special case where $\refplant{k}{k-1}$ is a maximally diffusive probability and $\refpolicy{k}{k-1}$ is uniform, as shown in~\cite{shafiei2025robustdecisionmakingfreeenergy}, the relaxed problem formulation yields a cost functional that aligns with the maximum diffusion objective in the MaxDiff framework~\cite{Berrueta2024}. Therefore, in this special case,  solving Problem~\ref{prob:main} yields a distributionally robust maximum diffusion policy. This connects our problem statement to the MaxDiff framework and, through it, to MaxEnt. In fact, MaxDiff is known to generalize MaxEnt policies and inherit its desirable robustness properties~\cite{eysenbach2022maximum}. Yet, to compute policies that robustly maximize a reward, MaxEnt must be used with a different/pessimistic objective. Moreover, the ambiguity set over which the policy is robust is not explicit but  emerges from the entropy of the optimal policy; explicit bounds for the ambiguity set are only available for discrete variables. In Problem~\ref{prob:main}, the ambiguity set is instead explicitly defined in the problem formulation, as a control requirement. Finally, when $\statecost{k}$ is a cost-to-go, the objective in Problem~\ref{prob:main} yields  the expected free energy minimized to determine actions in the context of active inference from computational neuroscience~\cite{FRISTON2009293} and this relates the formulation to other popular frameworks, such as risk-sensitive control, KL control, linear quadratic Gaussian regulator and control-as-inference.
In this sense, Problem~\ref{prob:main} offers a generalization of active inference~\cite{FRISTON2009293} accounting for ambiguity in the generative model. 

\section{Main Results}\label{sec:results}
First, we obtain a reformulation of Problem~\ref{prob:main} and characterize its key properties. This yields a resolution strategy where, first the free energy is maximized across the ambiguity set and then minimized in the policy space. Solving the maximization problem requires solving an infinite dimensional non-convex problem and tailoring a Frank-Wolfe algorithm to solve this problem. For brevity,  results are stated for continuous variables and we give a sketch of the proofs. Complete proofs will be presented elsewhere.

\subsection{Reformulating Problem~\ref{prob:main}}
The chain rule for the KL divergence (Lemma \ref{lem:chain_rule}) implies that the KL divergence term in the cost functional of~\eqref{eqn:main_problem} can be written as $\DKL{\policy{k}{k-1}}{\refpolicy{k}{k-1}} + \E_{\policy{k}{k-1}}\left[\DKL{\plant{k}{k-1}}{\refplant{k}{k-1}}\right]$. Thus, Problem~\ref{prob:main} can be reformulated as
\begin{multline}\label{eq:maingeneralproblem}
\underset{\shortpolicy{k}{k-1}\in\sD}{\min}   \left\{ \DKL{\shortpolicy{k}{k-1}}{\shortrefpolicy{k}{k-1}} + \E_{\shortjointxu{k}{k-1}}\left[l(\bv{x}_{k},\bv{u}_k)\right] \right. \\ \left.
 +  \underset{\shortrefplant{k}{k-1}\in\ball{k}{k-1}}{\max} \E_{\shortpolicy{k}{k-1}}\left[\DKL{\shortplant{k}{k-1}}{\shortrefplant{k}{k-1}}\right] \right\}.
\end{multline}
%\begin{eqnarray}\nonumber
%\underset{\shortpolicy{k}{k-1}\in\sD}{\min} & \DKL{\shortpolicy{k}{k-1}}{\shortrefpolicy{k}{k-1}} + \E_{\shortjointxu{k}{k-1}}\left[l(\bv{x}_{k},\bv{u}_k)\right] \\
%& +  \underset{\shortrefplant{k}{k-1}\in\ball{k}{k-1}}{\max} \E_{\shortpolicy{k}{k-1}}\left[\DKL{\shortplant{k}{k-1}}{\shortrefplant{k}{k-1}}\right].
%\end{eqnarray}
%l(\bv{x}_{k},\bv{u}_k) = \statecost{k} + \actioncost{k}$
The outer minimization is over policies $\shortpolicy{k}{k-1}$, while the inner maximization is over generative models $\shortrefplant{k}{k-1}$.
The next result shows that the feasibility set of the inner maximization is compact. 
\begin{lemma}
\label{lem:compactness} 
pick any $\shortnominalplant{k}{k-1}\in\sD$. The ambiguity set $\mathcal{B}_{\eta}(\shortnominalplant{k}{k-1})$ is compact in $L^1(\Omega)$.
\end{lemma}
\begin{proofsketch}
assume that $\shortnominalplant{k}{k-1} \in L^\infty(\Omega)$. By definition of the ambiguity set, for any $C > 1$, $\int_{\{q^{(x)}_{k\mid k-1}> C\}} q^{(x)}_{k\mid k-1}\le \frac{\eta_k(x_{k-1},u_k) + \|\bar q^{(x)}_{k\mid k-1}\|_{\infty}}
{\ln C}$, so the family of densities in $\mathcal{B}_{\eta}(\shortnominalplant{k}{k-1})$ is uniformly integrable.
By \cite[Theorem 4.7.20]{bogachev2007measure}, any uniformly integrable subset of $L^1(\Omega)$ whose defining constraints are closed is relatively compact in the $L^1(\Omega)$ weak topology.
Since the KL-divergence constraint and
non-negativity constraint are closed, the ambiguity set is also weakly closed and therefore compact.
\end{proofsketch}

Building on this result, we next show that the inner maximization problem admits an optimal solution and that the maximization and expectation in~\eqref{eq:maingeneralproblem} can be swapped.
\begin{theorem}\label{theo-swap}
pick any $\shortpolicy{k}{k-1}\in\sD$. The inner maximization problem in~\eqref{eq:maingeneralproblem}
admits an optimal solution $q^{(x),\ast}_{k|k-1}$, and
\begin{equation}
\begin{aligned}\label{eq:swapped}
&\underset{\shortrefplant{k}{k-1}\in \ball{k}{k-1}}{\max}\E_{\shortpolicy{k}{k-1}}\left[\DKL{\shortplant{k}{k-1}}{\shortrefplant{k}{k-1}}\right]\\
&\quad =\E_{\shortpolicy{k}{k-1}}\left[\underset{\shortrefplant{k}{k-1}\in \ball{k}{k-1}}{\max}\DKL{\shortplant{k}{k-1}}{\shortrefplant{k}{k-1}}\right].
\end{aligned}
\end{equation}
\end{theorem}

\begin{proofsketch}
fix $p^{(u)}_{k\mid k-1}$ and let $\bv{Z}:=(\bv{X}_{k-1},\bv{U}_k)$. By Lemma \ref{lem:compactness}, for each realization $\bv{z}$, the set $\Gamma(\bv{z}):=B_{\eta}\!\bigl(\bar q^{(x)}_{k\mid k-1}(\cdot\mid \bv{z})\bigr)$ is compact in $L^1(\Omega)$.  
Define the continuous function $f(\bv{z},\shortrefplant{k}{k-1})=D_{\mathrm{KL}}\!\bigl(p^{(x)}_{k\mid k-1}(\cdot\mid \bv{z})\,\|\,\shortrefplant{k}{k-1}\bigr)$.
Since $f(\bv{z},\cdot)$ is continuous on $\Gamma(\bv{z})$, following Weierstrass Theorem \cite[Corollary 2.35]{Aliprantis2006}, a maximizer for the problem on the right-hand side of~\eqref{eq:swapped} exists.
We denote the optimizer by $q^{(x),*}_{k\mid k-1}:=q^{(x),*}_{k}(\bv{x}_k\mid\bv{x}_{k-1},\bv{u}_k)\in\Gamma(\bv{z})$. Moreover, it can be shown that $\Gamma(\cdot)$ has measurable graph and $f$ is Carathéodory. Thus, a measurable selector $q^*(\bv{Z})$ exists \cite[Theorem 18.19]{Aliprantis2006}. This implies that (details omitted for brevity) $\max_{q^{(x)}_{k\mid k-1}\in\Gamma(Z)}
\mathbb{E}_{p^{(u)}_{k\mid k-1}}\!\bigl[f(Z,q^{(x)}_{k\mid k-1})\bigr]
=
\mathbb{E}_{p^{(u)}_{k\mid k-1}}\!\Bigl[\max_{q^{(x)}_{k\mid k-1}\in\Gamma(Z)} f(Z,q^{(x)}_{k\mid k-1})\Bigr]$, yielding the desired conclusion.
\end{proofsketch}

The above results imply that Problem~\ref{prob:main} can be solved by first finding $q^{(x),*}_{k\mid k-1}$, i.e.,  the optimal solution of the problem in the right hand side of~\eqref{eq:swapped} and then by minimizing in the policy space the functional
\begin{multline}\label{eq:functional_pu}
\DKL{\shortpolicy{k}{k-1}}{\shortrefpolicy{k}{k-1}} + \E_{\shortjointxu{k}{k-1}}\left[\statecost{k} + \actioncost{k}\right] \\
 +  \E_{\shortpolicy{k}{k-1}}\left[\DKL{\shortplant{k}{k-1}}{q^{(x),*}_{k\mid k-1}}\right].   
\end{multline}
%\begin{eqnarray}\nonumber
%& \DKL{\shortpolicy{k}{k-1}}{\shortrefpolicy{k}{k-1}} + \E_{\shortjointxu{k}{k-1}}\left[\statecost{k} + \actioncost{k}\right] \\
%& +  \E_{\shortpolicy{k}{k-1}}\left[\DKL{\shortplant{k}{k-1}}{q^{(x)*}_{k\mid k-1}}\right].
%\end{eqnarray}
The next result gives the optimal solution to this last problem.
\begin{lemma}\label{lem:optimal_policy}
consider the problem $\underset{\policy{k}{k-1}\in\sD}{\min} \mathcal{J}(\policy{k}{k-1})$, with $\mathcal{J}(\cdot)$ given in~\eqref{eq:functional_pu}. The optimal solution $\optimalpolicy{k}{k-1}$ is
\begin{equation}\label{eqn:optimal_policy}
\frac{\shortrefpolicy{k}{k-1}\exp\left(-\mathcal{T} (\bv{X}_k,\bv{u}_k)-\DKL{\shortplant{k}{k-1}}{q^{(x),*}_{k\mid k-1}}\right)}{Z},
\end{equation}
where $\mathcal{T} (\bv{X}_k,\bv{u}_k) \coloneqq\E_{\shortplant{k}{k-1}}\left[\statecostexpectation{k}\right]+\actioncost{k}$, and $Z$ is the normalizer.
\end{lemma}

\begin{proofsketch} follows from Lemma~\ref{lem:free_energy_optimization} and is omitted here for brevity.
%
%$\E_{\shortpolicy{k}{k-1}}\left[\E_{\shortplant{k}{k-1}}\left[\statecostexpectation{k}\right] + \actioncost{k} +\DKL{\shortplant{k}{k-1}}{q^{(x)*}_{k\mid k-1}} \right]$.
%$   \E_{\shortpolicy{k}{k-1}}\left[\E_{\shortplant{k}{k-1}}\left[\statecostexpectation{k}\right] + \actioncost{k}+\DKL{\shortplant{k}{k-1}}{q^{(x)*}_{k\mid k-1}} \right]$.
\end{proofsketch}

\subsection{Tackling The Maximization Problem}
Lemma~\ref{lem:optimal_policy} implies that, to compute the optimal policy, we need to find the optimal value for the problem
\begin{equation}\label{eqn:max_DKL}
    \underset{\shortrefplant{k}{k-1}\in\ball{k}{k-1}}{\max}\DKL{\shortplant{k}{k-1}}{\shortrefplant{k}{k-1}}.
\end{equation}
Inspired by~\cite{shafiei2025robustdecisionmakingfreeenergy,hu2013kullback}, to tackle this problem we introduce a change of variables, using the likelihood ratio (or Radon-Nikodym derivative). This is formalized with the following
\begin{definition}\label{def-ratios}
the likelihood ratio generated by $\shortnominalplant{k}{k-1}$ is $r_{k|k-1}\coloneqq\frac{q^{(x)}_{k|k-1}}{\bar{q}^{(x)}_{k|k-1}}$. The set of ratios is
$\sR\coloneqq\{\ratio{k}{k-1}:\E_{\shortnominalplant{k}{k-1}}\left[\ratio{k}{k-1}\right]=1, \ratio{k}{k-1}\ge0, \text{ a.s.}\}$.
\end{definition}
\begin{assumption}\label{ass:bounds-r}given $\bv{x}_{k-1}$ and $\bv{u}_k$, $\mathcal{R} \subset L^1(\Omega)$ and there exist constants $\delta_0, \delta_1 \in \mathbb{R}_+$ independent of $r_{k|k-1}\in \mathcal{R}$ such that $0 < \delta_0 \le r_{k|k-1} \le \delta_1 < \infty$.
\end{assumption}

\begin{remark}
Assumption~\ref{ass:bounds-r} is used to prove convergence of the algorithm we propose to tackle the maximization problem~\eqref{eqn:max_DKL}. This assumption is satisfied whenever, e.g., $\shortrefplant{k}{k-1}$ and $\shortnominalplant{k}{k-1}$ are Gaussians.
\end{remark}

With Definition~\ref{def-ratios}, the problem in~\eqref{eqn:max_DKL} becomes 
\begin{align}\label{eq:mainproblem-1}
&\E_{\shortnominalplant{k}{k-1}}\left[\frac{\shortplant{k}{k-1}}{\shortnominalplant{k}{k-1}}\ln\frac{\shortplant{k}{k-1}}{\shortnominalplant{k}{k-1}}\right]\nonumber - \underset{\ratio{k}{k-1}}{\min}
\E_{\shortnominalplant{k}{k-1}}\left[\frac{\shortplant{k}{k-1}}{\shortnominalplant{k}{k-1}}\ln \ratio{k}{k-1}\right],\\
&\text{s.t.} \quad \E_{\shortnominalplant{k}{k-1}}[\ratio{k}{k-1}\ln \ratio{k}{k-1}]\le\radius{k}{k-1}.
\end{align}
The first term in~\eqref{eq:mainproblem-1} does not depend on $\ratio{k}{k-1}$  and therefore the optimal ratio $\optimalratio{k}{k-1}$ can be found by solving the minimization in~\eqref{eq:mainproblem-1}. This is a minimization of a concave function, which we tackle by proposing a reformulation within the DC programming framework. In fact, we have:
\begin{equation}\label{eqn:difference_reformulation}
    \ln \ratio{k}{k-1}=\ratio{k}{k-1}\ln \ratio{k}{k-1}-(\ratio{k}{k-1}\ln \ratio{k}{k-1}-\ln \ratio{k}{k-1}).
\end{equation} 
In~\eqref{eqn:difference_reformulation} we expressed the logarithm in the cost in~\eqref{eq:mainproblem-1} as the difference of two convex functions. This leads to the following DC reformulation of the problem in~\eqref{eq:mainproblem-1}
\begin{problem}\label{problem_main}
%assume that Assumption \ref{ass:bounds-r} holds.
find $\optimalratio{k}{k-1} \in \sR$ that minimizes
\begin{equation}\label{eq:mainproblem}
    \phi_k(\ratio{k}{k-1}) \coloneqq
g_k(\ratio{k}{k-1})-h_k(\ratio{k}{k-1}),
\end{equation}
with
\begin{gather}\label{func_gk}
g_k(\ratio{k}{k-1})\coloneqq\E_{\shortnominalplant{k}{k-1}}\left[\frac{\shortplant{k}{k-1}}{\shortnominalplant{k}{k-1}}\ratio{k}{k-1}\ln \ratio{k}{k-1} \right],\text{ and}\\
h_k(\ratio{k}{k-1}) \coloneqq\E_{\shortnominalplant{k}{k-1}}\left[\frac{\shortplant{k}{k-1}}{\shortnominalplant{k}{k-1}} \kappa(\ratio{k}{k-1})  \right], \label{func_hk}
%\text{ where} \\ \nonumber\kappa(\ratio{k}{k-1}) \coloneqq\ratio{k}{k-1}\ln \ratio{k}{k-1}-\ln \ratio{k}{k-1}.
\end{gather}
where $\kappa(\ratio{k}{k-1}) \coloneqq\ratio{k}{k-1}\ln \ratio{k}{k-1}-\ln \ratio{k}{k-1}$,  subject to $\ratio{k}{k-1} \in \sC$, with $\sC \subseteq \sR$ given by
%{\color{blue}---the constraint*** unclear -- the ratio should belong to the set C. Here, you can just also write that it is a subset of R}
\begin{equation}\label{eq:mainproblem2}
\mathcal{C} \coloneqq \E_{\shortnominalplant{k}{k-1}}[\ratio{k}{k-1}\ln \ratio{k}{k-1}]\le\radius{k}{k-1}.
\end{equation}
\end{problem}
Next, we propose an algorithm to solve the DC Problem~\ref{problem_main}. Properties of the algorithm are characterized in Section~\ref{sec:properties}.
\begin{remark}
the definition of ambiguity set ensures that the optimal value of the maximization problem is bounded.    
\end{remark}

\begin{lemma}\label{lemma:C-compact}
assume that Assumption \ref{ass:bounds-r} holds.
Then $\sC$ in Problem \ref{problem_main} is compact (and convex and closed) in $L^1(\Omega)$.
\end{lemma}

\begin{proofsketch}
from Assumption \ref{ass:bounds-r}, the bounds $\delta_0 \le r_{k\mid k-1} \le \delta_1$ imply that $\mathcal{R}$ is
bounded in $L^1(\Omega)$, and convexity follows from the bounds being
preserved under convex combinations.
Since $\mathcal{R}$ is bounded and closed in
$L^1(\Omega)$, it is compact, which proves the claim.
Assumption \ref{ass:bounds-r} also implies that   the feasible set $\sC$ is obtained as the intersection of closed sets. This yields the result.
\end{proofsketch}

\subsection{Tackling Problem~\ref{problem_main}}\label{sec:Alg}

The algorithm we propose to tackle Problem~\ref{problem_main} is a Frank-Wolfe algorithm~\cite{frank1956} for DC programs adapted from~~\cite{millan2023frank} to a constrained infinite dimensional setting.

Before describing the key steps of the algorithm -- given in Algorithm \ref{alg:svrfi} -- we first show that the functions $g_k(\cdot)$ and $h_k(\cdot)$ in~\eqref{eq:mainproblem} are Lipschitz continuous over the feasibility set. Lipschitz continuity of these function implies that the objective in~\eqref{eq:mainproblem} is also Lipschitz. These properties are crucial to characterize convergence of the algorithm. 

\begin{theorem}\label{theo:Lipschitz_phi}
let $\rho_{k|k-1}^{(x)} \coloneqq \frac{\shortplant{k}{k-1}}{\shortnominalplant{k}{k-1}}$ and assume that: (i) $\rho \in L^\infty(\Omega)$, with $\|\rho\|_\infty \le M$; (ii)  Assumption \ref{ass:bounds-r} holds.
Then: (i) $g_k(r_{k|k-1})$ and $h_k(r_{k|k-1})$ in \eqref{func_gk} and \eqref{func_hk} are Lipschitz continuous on $\mathcal{R}$ with Lipschitz constants being $L_{g_k} \le M (|\ln \delta_0| + 1)$ and $L_{h_k} \le M  \max \{ |\ln \delta_0|, |\ln \delta_1| \}$, respectively; (ii) $\phi_k(\ratio{k}{k-1}) \coloneqq g_k (\ratio{k}{k-1})- h_k(\ratio{k}{k-1})$ is  Lipschitz continuous on $\mathcal{R}$, with constant $L_{\phi_k} \le L_{g_k} + L_{h_k}$.
\end{theorem}

\begin{proofsketch}
by Assumption \ref{ass:bounds-r}, the likelihood ratio satisfies
$r_{k\mid k-1} \in [\delta_0,\delta_1]$. On this interval, the scalar functions $\psi_1(r)=r\ln r$ and $\psi_2(r)=r\ln r-\ln r$ have bounded derivatives and hence Lipschitz. This also implies that their difference is Lipschitz. The specific constants in the statement are obtained by leveraging the definition of the functions and the fact that $\|\rho^{(x)}_{k\mid k-1}\|_\infty \le M$.
\end{proofsketch}

\begin{algorithm}
  \caption{Frank-Wolfe algorithm for Problem~\ref{problem_main}}
  \label{alg:svrfi}
  \begin{algorithmic}[1]
      \STATE {\bf Inputs:} $g_k(\cdot)$, $h_k(\cdot)$, $\shortplant{k}{k-1}$, $\shortnominalplant{k}{k-1}$ and $\radius{k}{k-1}$.
    \STATE {\bf Initialize:} $n \gets 0$; $L_0 \gets L_{g_k}  + L_{h_k}$ (Theorem~\ref{theo:Lipschitz_phi}); select $r^n_{k|k-1}\in\mathcal{C}$.
  %  \FOR{$n=1$ \TO $T$}
        \STATE Compute $v^n_k \gets \nabla h_k(r^n_{k|k-1})$, and $p_{k|k-1}^n$ solving 
\begin{gather} \label{eq_p^n_algo}
    p^n_{k|k-1} =\mathop{\argmin}\limits_{p_{k|k-1}\in\mathcal{C}} \left\langle {\nabla g_k(r^n_{k|k-1})-v_k^n},p_{k|k-1}-r^n_{k|k-1}\right\rangle.
\end{gather}
      \STATE Compute  $ \omega_k(r^n_{k|k-1})$ as
      \begin{equation}
    \label{eq_omega}
    \omega_k(r^n_{k|k-1}) \gets \left\langle {\nabla g_k(r^n_{k|k-1})-v_k^n},p_{k|k-1}^n-r^n_{k|k-1}\right\rangle.
\end{equation}
      \WHILE{$\omega_k(r^n_{k|k-1}) \ne 0$}
        \STATE Compute $j$ as
        \begin{equation}\label{eq_j_alg}
        j \gets \min\{l\in\N : 2^lL_n\ge2L_0\},
        \end{equation}
        and the step size {$\lambda^j_k \in (0,1]$}, given by
        \begin{align}\label{eq_alg_Lambda}
\lambda^j_k = \min\left\{1,\frac{|\omega_k(r^n_{k|k-1})|}{2^jL_n\|p^n_{k|k-1}-r^n_{k|k-1}\|_{L^1}^2}\right\}.
\end{align}
        \IF{Inequality \eqref{eq_ineq_step4} holds}
          \STATE $j_n \gets j$
          \STATE $\lambda^n_k \gets \lambda_k^{j_n}$
          \STATE Compute $r^{n+1}_{k|k-1}$ and $L_{n+1}$ as
          \begin{align}
\label{eq:alg_r}
&r^{n+1}_{k|k-1}
= r^n_{k|k-1} + \lambda^n_k (p^n_{k|k-1} - r^n_{k|k-1}), \\
\label{eq:alg_L}
&L_{n+1}= 2^{j_n - 1} L_n.
\end{align}
\STATE Compute $v^{n+1}_k$ (line 3) and $p_{k|k-1}^{n+1}$ using~\eqref{eq_p^n_algo} with $r^{n+1}_{k|k-1}$ and $v^{n+1}_k$.
\STATE Compute $\omega_k(r^{n+1}_{k|k-1})$ using~\eqref{eq_omega} with $r^{n+1}_{k|k-1}$, $v^{n+1}_k$
          \STATE $n \gets n+1$
        \ELSE
        \WHILE{Inequality \eqref{eq_ineq_step4} does not hold}
          \STATE $j \gets j+1$
          \STATE Compute the step size $\lambda_k^j$ as in \eqref{eq_alg_Lambda}. 
        \ENDWHILE
        \STATE Repeat steps from lines $8$ to $13$.  
        \ENDIF
      \ENDWHILE
      \STATE \textbf{return} $r_{k|k-1}^\ast \gets r^n_{k|k-1}$.
 %   \ENDFOR
  \end{algorithmic}
\end{algorithm}

Now, we describe Algorithm~\ref{alg:svrfi}. In Section~\ref{sec:properties} we show that the algorithm returns a sequence of ratios, say $r_{k|k-1}^n$, with $n$ indexing the algorithm iterations, converging to the optimal solution of Problem \ref{problem_main}. Sequence elements are obtained by evaluating a linear {\em oracle} -- returning the optimal solution of a linear problem -- and subsequently updating the step size through a backtracking rule.

%\textcolor{blue}{Where is L determined? It is needed early in the algorithm. Comment on the inputs, e.g., $L_0$.}
The algorithm takes as input $g_k(\cdot)$ and $h_k(\cdot)$, together with $\shortplant{k}{k-1}$, $\shortnominalplant{k}{k-1}$ and $\radius{k}{k-1}$.
The procedure is initialized by selecting the initial iterate $r^0_{k|k-1} \in \sC$ and the constant $L_0 = L_{g_k}  + L_{h_k}$ given by Theorem~\ref{theo:Lipschitz_phi}.
Also, $j$ denotes the backtracking index used to adaptively adjust the estimate Lipschitz constant $L_n$, $\lambda_k^j$ is the candidate step size, and $\lambda_k^n$ is the accepted step size.
In the spirit of the standard Frank-Wolfe, Algorithm \ref{alg:svrfi}: (i) solves the linear subproblem  to obtain the search direction (line 3); (ii) computes a candidate step–size (line 6), the descent condition (line 7), and the final update of the iterate (line 10). 
In particular, the algorithm first (line 3) computes the Frank–Wolfe search direction by solving a linear minimization problem yielding $p_{k|k-1}^n$.
Essentially, this step picks the feasible point that minimizes the linear model around the current iterate.
%together with an index $j$
Then, $\omega_k(r^n_{k|k-1})$ is computed (defined in line 4) and the algorithm checks if $r^n_{k|k-1}$ satisfies the first-order optimality condition (line 5). If this condition is not met, an index $j$ and a candidate step size $\lambda_k^j$ are computed (line 6).
Then,  the following inequality is evaluated (line 7) 
\begin{align}\nonumber
&\phi_k\big(r^n_{k|k-1} + \lambda^j_k (p^n_{k|k-1} - r^n_{k|k-1})\big)
\le \phi_k(r^n_{k|k-1}) \\ \label{eq_ineq_step4}
&- |\omega_k(r^n_{k|k-1})| \lambda^j_k  + \frac{2^j L_n}{2} \|p^n_{k|k-1} - r^n_{k|k-1}\|_{L^1}^2 (\lambda^j_k)^2.
\end{align}
If the above quadratic majorization holds,  the algorithm updates its internal variables for the next iteration (lines 10 -- 13).
Otherwise (lines 15--19), the step size is rejected and new candidates are generated until \eqref{eq_ineq_step4} is met; once this descent condition holds, the internal variables are subsequently updated.
The algorithm outputs its final iterate $r^n_{k|k-1}$.

{\begin{remark}
the backtracking index
$j$ is used to adaptively adjust the estimate of the curvature of the objective function.
By increasing $j$, the algorithm ensures that the quadratic upper bound in \eqref{eq_ineq_step4} is satisfied. As will be formalized in Section \ref{sec:properties}, this mechanism guarantees that the backtracking procedure terminates and that each accepted step yields a decrease in the objective function.
\end{remark}
}

\subsection{Properties of Algorithm \ref{alg:svrfi}}\label{sec:properties}
%\textcolor{blue}{It is uncleear how the quantities in the algorithm are related to the theorems here}

%We show that the algorithm is well-posed and that its iterates behave in a predictable manner.
%We rely on previously established results to guarantee that each step of the algorithm is feasible, that the objective decreases along accepted iterations, and that the duality gap provides a meaningful measure of progress toward optimality. Together, these properties confirm that Algorithm \ref{alg:svrfi} is consistent with Problem \ref{problem_main} and converges to a stationary solution.

To begin with, the next lemma guarantees that the gradients used in the linearization steps (lines 3 and 4) are Lipschitz on the feasible domain $\sC$.
Therefore, the linear oracle always has a solution, and the quadratic majorization~\eqref{eq_ineq_step4} used in the backtracking step is valid.

\begin{lemma}\label{theo_nabla_Lipschitz}
let $\rho_{k|k-1}^{(x)} \coloneqq \frac{\shortplant{k}{k-1}}{\shortnominalplant{k}{k-1}}$, $\sR$ as in Definition \ref{def-ratios} and $\mathcal{C}$ given by \eqref{eq:mainproblem2} and assume that: (i) $\rho \in L^\infty(\Omega)$ with $\|\rho\|_\infty < \infty$; (ii)  Assumption \ref{ass:bounds-r} holds.
Then, $\nabla g_k(\ratio{k}{k-1})$, $\nabla h_k(\ratio{k}{k-1})$ and $\nabla \phi_k(\ratio{k}{k-1})$ are Lipschitz on $\sR$ (hence also on $\sC \subseteq \sR$).
\end{lemma}

\begin{proofsketch}
the functions involved are $g_k(\cdot)$ and $h_k(\cdot)$ as in \eqref{func_gk} and \eqref{func_hk}, respectively, with gradients $\nabla g_k(\ratio{k}{k-1})=\rho_{k|k-1}^{(x)}(\ln \ratio{k}{k-1} +1 )$ and $\nabla h_k(\ratio{k}{k-1})=\rho_{k|k-1}^{(x)} \left(\ln \ratio{k}{k-1} +1 +\frac{1}{\ratio{k}{k-1}} \right)$.
Define the functions $\psi_1(\ratio{k}{k-1})=\ratio{k}{k-1} \ln \ratio{k}{k-1}$ and $\psi_2(\ratio{k}{k-1})=\ratio{k}{k-1} \ln \ratio{k}{k-1}-\ln \ratio{k}{k-1}$.
Both are $C^2$ on $(0,\infty)$, and their second derivatives are $\psi_1''(\ratio{k}{k-1})=\frac{1}{\ratio{k}{k-1}}$ and $\psi_2''(\ratio{k}{k-1})=\frac{1}{\ratio{k}{k-1}}+\frac{1}{\ratio{k}{k-1}^2}$.
By Assumption \ref{ass:bounds-r} every $\ratio{k}{k-1} \in \sR$ satisfies $\delta_0 \le r_{k|k-1} \le \delta_1$, a.s.
Hence both second derivatives $\psi_1''$ and $\psi_2''$ are uniformly bounded on the interval $[\delta_0, \delta_1]$.
Therefore, $\psi_1''$ and $\psi_2''$ are Lipschitz on $[\delta_0, \delta_1]$.
Since $\rho \in L^\infty(\Omega)$, multiplying a Lipschitz function by a bounded function preserves Lipschitz continuity.
Thus $\nabla g_k(r_{k|k-1})=\rho \psi_1'(r_{k|k-1})$ and $\nabla h_k(r_{k|k-1})=\rho \psi_2'(r_{k|k-1})$ are Lipschitz on $\sR$.
In addition, $\nabla \phi (r_{k|k-1}) = \nabla g_k(r_{k|k-1}) - \nabla h_k(r_{k|k-1})$, so Lipschitz continuity of $\nabla \phi (r_{k|k-1})$ follows.
Since $\sC \subseteq \sR$ the same Lipschitz property holds automatically on $\sC$.
\end{proofsketch}

Next, we show that if a given functional $\psi(\cdot)$ has an
$L$-Lipschitz gradient over the feasible region, then an update of the form $f(\cdot) + t g(\cdot)$ yields an explicit upper bound in
terms of the directional derivative and a quadratic curvature term. This inequality plays a central role in establishing the convergence of Algorithm~\ref{alg:svrfi} in Theorem \ref{theo_alg1}.

\begin{theorem}\label{theo_nabla_psi}
let $\psi:L^1(\Omega)\to\R$ be continuously differentiable with $L$-Lipschitz gradient on $\Omega$. Assume that $f\in\Omega$, $g\in L^1$, $t\in[0,1]$ and $f+tg\in\Omega$.
Then, $\psi(f+tg)\le\psi(f)+t \left\langle \nabla\psi(f),g\right\rangle+\frac{L}{2}\|g\|_{L^1}^2t^2$.
\end{theorem}

\begin{proofsketch}
set $h(t):=\psi(f+t g)$, $t\in[0,1]$.
By chain rule,
$h'(t)=\langle \nabla\psi(f+t g),\, g\rangle$. 
Then $\psi(f+t g)-\psi(f)=\int_0^t \langle \nabla\psi(f+s g),\, g\rangle\, ds
= t\langle \nabla\psi(f),g\rangle
+ \int_0^t \langle \nabla\psi(f+s g)-\nabla\psi(f),\, g\rangle\, ds$.
Using Cauchy--Schwarz and the $L$-Lipschitz property of $\nabla\psi$, $\|\nabla\psi(f+s g)-\nabla\psi(f)\|_{L^1} \le Ls\|g\|_{L^1}$, so $\psi(f+t g)\le \psi(f)+t\langle \nabla\psi(f),g\rangle+\tfrac{L}{2}\|g\|_{L^1}^2 t^2$, and the proof is complete.
\end{proofsketch}

%Theorem \ref{decreasing} implies that the backtracking  in
%Algorithm~\ref{alg:svrfi} always \textcolor{blue}{accepts a step} if -- as enforced in Line $6$ -- $2^j L_n \ge 2L_0$. In turn, this ensures that the update based on $\lambda_k^{\,j}$ is well defined.
%Moreover, the above result also shows that each accepted step yields a strict decrease of
%$\phi_k$, which implies monotonic progress of the sequence
%$\{\phi_k(r^n_{k\mid k-1})\}$.

Inequality \eqref{eq_ineq_step4} ensures that the objective value decreases by at least a quadratic amount in the step size.
Therefore, combined with the convex combination update \eqref{eq:alg_r} and \eqref{eq:alg_L}, the sequence
$\{\phi_k(r^n_{k\mid k-1})\}$ is strictly monotone and upper bounded, which implies convergence.
This property is formalized with the next result.

\begin{remark}
from Theorem \ref{theo:Lipschitz_phi} and Lemma \ref{theo_nabla_Lipschitz}, both $g_k(\cdot)$ and $\nabla g_k(\cdot)$ are Lipschitz on $\sC$, with  $g_k(\cdot)$ as in \eqref{func_gk}. 
Hence, Theorem \ref{theo_nabla_psi} holds for  $g_k(\cdot)$ giving us 
\begin{align}\nonumber
 &g_k(r^n_{k|k-1}+\lambda^j_k(p^n_{k|k-1}-r^n_{k|k-1})) \\ \nonumber
 &\le g_k(r_{k|k-1}^n) +\left\langle \nabla g_k(r_{k|k-1}^n),p^n_{k|k-1}-r^n_{k|k-1}\right\rangle\lambda^j_k\\ \label{ineq_gk}
 &+\frac{L_{g_k}}{2}\|p^n_{k|k-1}-r^n_{k|k-1}\|_{L^1}^2\lambda_k^j,
\end{align}\nonumber
with $j \coloneqq \min\{l\in\N : 2^lL_n\ge2L_0\}$ and
\begin{align}
\lambda^j_k & \coloneqq\mathop{\argmin}\limits_{\lambda\in(0,1]}\left\{-|\omega_k(r^n_{k|k-1})|\lambda+\frac{2^jL_n}{2}\|p^n_{k|k-1}-r^n_{k|k-1}\|_{L^1}^2\lambda^2\right\}\\
&=\min\left\{1,\frac{|\omega_k(r^n_{k|k-1})|}{2^jL_n\|p^n_{k|k-1}-r^n_{k|k-1}\|_{L^1}^2}\right\}.
\end{align}
As discussed in~\cite{millan2023frank}, the index $j$ defined as above automatically adjusts the curvature estimate via backtracking, while $\lambda^j_k$ provides a safe step that guarantees descent of the objective in the Frank–Wolfe method for DC problems.
Inequality \eqref{ineq_gk} will be essential in the proof of Theorem \ref{decreasing}.
\end{remark}

%{\color{blue}for the next theorem $L_{g_k}=2^jL_n$}.

\begin{theorem}\label{decreasing}
the following assertions are true:
\begin{enumerate}
    \item If $2^jL_n\ge 2L_0$, with $L_n$ (line 10) standing for an estimate of the Lipschitz constant of $\phi(\cdot)$, for some $j \in \N$, then inequality \eqref{eq_ineq_step4} holds;
    \item Let $j_n$ be the smallest non-negative integer satisfying  $2^{j_n}L_n\ge 2L_0$. Then, \eqref{eq_ineq_step4} also holds for $j=j_n$ and $\lambda_k^n=\lambda_k^{j_n}$ \eqref{eq_alg_Lambda}. Furthermore,
    \begin{align}
    \phi_k(r^{n+1}_{k|k-1})\le\phi_k(r^n_{k|k-1})-\frac{1}{2}\omega_k(r^n_{k|k-1})\lambda^n_k.
\end{align}
\end{enumerate}
\end{theorem}

\begin{proofsketch}
1) from \eqref{eq_p^n_algo} the search direction satisfies
$\omega_k(r^n_{k\mid k-1})
\le \langle \nabla \phi_k(r^n_{k\mid k-1}),\, p^n_{k\mid k-1} - r^n_{k\mid k-1} \rangle \le 0$.
Using the Lipschitz model implied by $2^j L_n \ge 2L_0$, the convex
combination update $L_{g_k}=2^jL_n$, in \eqref{ineq_gk}, and the convexity of function of $h_k(\cdot)$, \eqref{eq_ineq_step4} holds.
2) with $j=j_n$ and $\lambda_k^n=\lambda_k^{j_n}$, then \eqref{eq_ineq_step4} is satisfied directly from 1).
In addition, from \eqref{eq_ineq_step4} together with the updates of $\ratio{k}{k-1}^{n+1}$ and $L_{n+1}$ described in line 10 of Algorithm \ref{alg:svrfi}, we get
$\phi_k(r^{n+1}_{k\mid k-1})
\le \phi_k(r^n_{k\mid k-1})
- \tfrac{1}{2}\,\omega_k(r^n_{k\mid k-1})\,\lambda_k^{\,n}$, showing that each accepted step decreases $\phi_k$ and that the backtracking terminates.
\end{proofsketch}

Algorithm~\ref{alg:svrfi} relies on getting $p^n_{k|k-1} \in \sC$ by solving the linear problem in~\eqref{eq_p^n_algo}. Since the feasibility domain is compact (Lemma~\ref{lemma:C-compact}), the problem  in~\eqref{eq_p^n_algo} always admits a solution, and every iterate $r^{n+1}_{k|k-1}
= r^n_{k|k-1} + \lambda^n_k (p^n_{k|k-1} - r^n_{k|k-1})$, with $\lambda^n_k \in (0,1]$, generated by the algorithm remains feasible. This follows from the fact that $p^n_{k|k-1}$ is obtained by solving the linear problem on line 3 of the algorithm and $\mathcal{C}$ is convex.

We note that Algorithm~\ref{alg:svrfi} guarantees non–degenerate progress at each iteration. In particular, the adaptive Lipschitz estimate remains uniformly bounded and the step size never vanishes until optimality conditions are met.
These properties follow from the structure of the backtracking rule and the compactness of the feasible set. Since their detailed proof is not essential for the development of the main results, we omit it and will present the full derivation elsewhere.

%The next lemma guarantees that the algorithm keeps making progress: the Lipschitz estimate stays bounded and the step size never vanishes until optimality conditions are met.

%\begin{lemma}\label{lemma:algo1}
%Let $j$ be defined as in \eqref{eq_j_alg}.
%Then, each element of the sequence $\{L_n\}_{n\in N}$ in \eqref{eq:alg_L} satisfies $L_0\le L_n\le L+L_0$, with $L=2^jL_n$.
%Furthermore, the step size sequence $
%\{\lambda^n_k\}_{k\in\N}$ satisfies  
%\begin{align*}
%\lambda^n_k\ge\min\left\{1,\frac{|\omega_k(r^n_{k|k-1})|}{2^{n-1}\alpha}\right\},\text{ where}
%\end{align*}
%$\alpha=2(L+L_0)({\rm %diam}\mathcal{C})^2>0$,  ${\rm %diam}\mathcal{C}=\sup_{r,s\in\mathcal{C}}\|r-s\|_{L^1}$.
%\end{lemma}
%\begin{proofsketch}
%from Algorithm \ref{alg:svrfi}, together with equation \eqref{eq:alg_L}, we have that $L_n=2^{j_{n-1}-1}L_{n-1}=2^{j_{n-1}-1} \dots 2^{j_0-1}L_0\ge L_0$   for all $n$.
%Then, the following inequality holds $2^nL_n\|p^n_{k|k-1}-r^n_{k|k-1}\|^2_{L^1} \le  2^n(L+L_0)(diam\mathcal{C})^2 =2^{n-1}\alpha$.
%Hence, 
%\begin{align*}
%\frac{|\omega_k(r^n_{k|k-1})|}{2^jL_n\|p^n_{k|k-1}-r^n_{k|k-1}\|_{L^1}^2}\ge\frac{|\omega_k(r^n_{k|k-1})|}{2^{n-1}\alpha}.
%\end{align*}
%Consequently, we obtain $\lambda^n_k\ge\min\Big\{1,\frac{|\omega_k(r^n_{k|k-1})|}{2^{n-1}\alpha}\Big\}$, as a direct result of \eqref{eq_alg_Lambda}.
%\end{proofsketch}}

Next, we establish the convergence of Algorithm \ref{alg:svrfi}. Without this result, we would only be guaranteed that the algorithm decreases the objective
(Theorem \ref{decreasing}), but not that it converges to a stationary point.

\begin{theorem}\label{theo_alg1}
let $\{r_{k|k-1}^n\}$ denote the sequence generated by Algorithm \ref{alg:svrfi}.
Then, the sequence $\{\phi_k(r^n_{k|k-1})\}$ is monotone decreasing, and every limit point of $\{r^n_{k|k-1}\}$ is a solution to the problem $\min\limits_{r_{k|k-1}\in\mathcal{C}}\phi_k(r_{k|k-1})$.
\end{theorem}

\begin{proofsketch}
from Theorem \ref{decreasing}, the sequence $\{\phi_k(r^{n}_{k\mid k-1})\}$ is
monotonically decreasing and upper bounded, hence convergent.
This property implies that $\omega_k(r^{n}_{k\mid k-1})$ converges to zero.
Since $\mathcal{R}$ is compact (Lemma \ref{lemma:C-compact}), the sequence $\{ r^{n}_{k\mid k-1} \}$ has limit points,
and any such limit point satisfies $\omega_k(r_{k\mid k-1}) = 0$,
which is the first-order optimality condition of Problem~\eqref{problem_main}. Thus, every accumulation
point of the sequence is a stationary solution of the problem.
\end{proofsketch}

In what follows, we show that at each iteration of Algorithm~\ref{alg:svrfi}, the step–size $\lambda_k^n$ remains strictly positive whenever $\omega_k(r^n_{k\mid k-1})$ is non-zero.
This prevents stagnation of the algorithm and ensures that it progresses towards the fulfillment of the first–order optimality condition.

\begin{lemma}\label{Lemma_decrease_phi}
consider the setting of Lemma \ref{theo_nabla_Lipschitz}.
Let $\phi_k$ be continuously differentiable with Lipschitz gradient and Lipschitz constant $L_{\nabla \phi_k}>0$ over a compact set $\mathcal{C}$.
Let $\{r_{k|k-1}\}$ denote the sequence generated by Algorithm \ref{alg:svrfi}.
Then:
\begin{align}\nonumber
&\phi_k(r^n_{k|k-1})-\phi_k(r^{n+1}_{k|k-1}) \\
&\ge\frac{\omega_k^2(r^n_{k|k-1})}{2\|r^n_{k|k-1}-p^n_{k|k-1}\|_{L^1}^2}\min\left\{\frac{1}{L_{\nabla \phi_k}},\frac{\|r^n_{k|k-1}-p^n_{k|k-1}\|_{L^1}}{\|\nabla\phi_k(r^n_{k|k-1})\|_{L^1}}\right\}
\end{align}
\end{lemma}

\begin{proofsketch}
using the Lipschitz continuity of $\nabla \phi_k(\ratio{k}{k-1})$, we have that $\phi_k(\ratio{k}{k-1}^n+\lambda(p^n_{k|k-1} - \ratio{k}{k-1}^n)) \leq \phi_k(\ratio{k}{k-1}^n)+\lambda \left\langle \nabla \phi_k(\ratio{k}{k-1}^n), p^n_{k|k-1} - \ratio{k}{k-1}^n\right\rangle + 2^{-1}L_{\nabla \phi_k} \lambda^2 \| p^n_{k|k-1} - \ratio{k}{k-1}^n \|^2$.
Since the linear minimization oracle ensures $\omega_k(\ratio{k}{k-1}^n)=\left\langle \nabla \phi_k(\ratio{k}{k-1}^n), p^n_{k|k-1} - \ratio{k}{k-1}^n\right\rangle$, substituting this direction into the inequality above yields a quadratic expression in $\lambda$.
Minimizing it over $\lambda \in (0,1]$ produces an explicit positive decrease in $\phi_k(\cdot)$.
\end{proofsketch}

%\textcolor{blue}{Need a paragraph at the beginning linking all the results, Need to link the results to the algorithm steps.}

\section{Numerical Example}

\begin{figure}[thpb]
\centering
    \includegraphics[width=\linewidth]
    {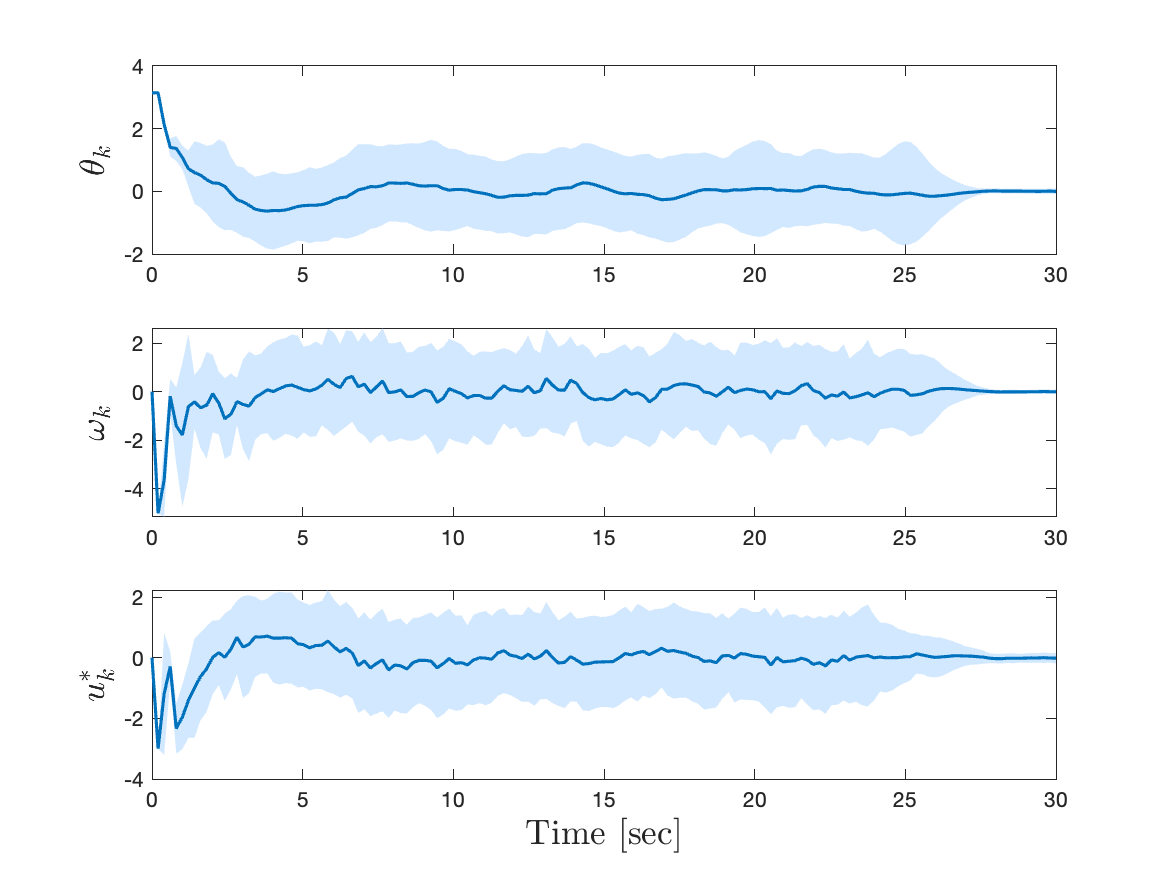}
    %{Figures/figtrajectories.eps}
    \caption{Simulation results. Panels show time evolution of $\theta_k$, $\omega_k$, and $u^\ast_k$ sampled from $\shortoptimalpolicy{k}{k-1}$. Solid lines represent the mean, and the shaded region represents the standard deviation obtained from $10^2$ simulations.}
    \label{fig:pendulum_results}
\end{figure}

As a representative example, to illustrate the effectiveness of our results, we consider the problem of stabilizing the unstable equilibrium  $[0,0]^\top$ of a pendulum, modeled by%
\begin{equation}
\begin{aligned}
\theta_{k} &= \theta_{k-1} + \omega_{k-1} dt + W_{\theta},\\
\omega_{k} &= \omega_{k-1} + \Big(\frac{g}{l}\sin(\theta_{k-1}) + \frac{u_k}{m l^2}\Big)dt + W_{\omega},
\end{aligned}
\label{eq:pendulum_dyn}
\end{equation}
where $\theta_k$ is the angular position, $\omega_k$ is the angular velocity, and $u_k$ is the applied torque. Also, $W_{\theta} \sim \mathcal{N}(0,0.05)$, $W_{\omega} \sim \mathcal{N}(0,0.1)$ so that $\shortplant{k}{k-1}$ is a Gaussian centered in~\eqref{eq:pendulum_dyn} with covariance matrix {$\Sigma_1 = \begin{bmatrix} 0.05 & 0 \\ 0 & 0.1\end{bmatrix}$}. In our experiments, the state is $\bv{X}_k = \begin{bmatrix}  \theta_k & \omega_k \end{bmatrix}^\top$ and $\mathcal{X} \coloneqq [-\pi,\pi] \times [-5,5]$ is discretized in a $20 \times 20$ grid. Also, $\mathcal{U} \coloneqq [-3,3]$ and this is discretized in a $50\times 1$ grid. The parameters of the pendulum are: $g= 9.81\;\mathrm{m/s^2}$, $l= 0.50\;\mathrm{m}$, $m = 0.50\;\mathrm{kg}$, $dt = 0.05 \;\mathrm{s}$.
Here, ambiguity arises in the generative model: we optimize over all $\shortrefplant{k}{k-1}$ inside $\ball{k}{k-1}$, with radius $\eta_k(\bv{x}_{k-1},u_k)$.
%\textcolor{blue}{say in few words where ambiguity arises. Bold Characters are again missing!}
%
The nominal generative model simply captures the desired configuration for the pendulum so that $\shortnominalplant{k}{k-1} = \mathcal{N}(0,\Sigma_2)$, with $\Sigma_2 = \begin{bmatrix} 0.1 & 0 \\ 0 & 0.01 \end{bmatrix}$ so that minimizing the KL divergence term in \eqref{eq:functional_pu} amounts at reducing the error from the target. The generative model is however ambiguous and this is captured via the ambiguity set $\ball{k}{k-1}$ with $\eta_k(\bv{X}_{k-1},u_k)= 0.3 + 0.02 \, \omega_{k-1}^2 + 0.01 \, {u}_k^2$.
%$\radius{k}{k-1} = 0.3 + 0.02 \, \omega_k^2 + 0.01 \, {u}_k^2$.
Intuitively, the constant term is a baseline and the other terms capture the fact that ambiguity is higher when the pendulum exhibits faster motion or requires larger torques. Finally, in the experiments we set $\statecost{k}=\theta_k^2+0.5\,\omega_k^2$, and $c_k^{(u)}(u_k)=0.1 {u}_k^2$. Given this set-up, the optimal policy is computed by solving Problem~\ref{prob:main}. This means, in accordance with Lemma~\ref{lem:optimal_policy} that the optimal policy is the softmax~\eqref{eqn:optimal_policy}. Optimal actions, $u_k^{\ast}$ are obtained by sampling from this policy and, in order to do so, $q^{(x),*}_{k\mid k-1}$ needs to be computed. This is the generative model solving the inner maximization problem in~\eqref{eqn:max_DKL}. To obtain $q^{(x),*}_{k\mid k-1}$, we leverage Algorithm~\ref{alg:svrfi}. In particular, using this algorithm, we first compute $r^\ast_{k \mid k-1}$ and then, by means of Definition~\ref{def-ratios}, we set $q^{(x),*}_{k\mid k-1} =\shortnominalplant{k}{k-1}r^\ast_{k \mid k-1}$\footnote{Our Matlab implementation is available upon request}. The results are illustrated in Figure~\ref{fig:pendulum_results}, which clearly shows that the actions sampled from the optimal policy in~\eqref{eqn:optimal_policy} stabilize the unstable equilibrium of the pendulum despite ambiguity. 

%{\color{red}
%We consider the control of a pendulum whose goal, in the forward problem,is to stabilize the unstable equilibrium at $\theta_d=0$ and $\omega_d=0$.
%The stage costs are
%\begin{gather}
%c(\bv{x}_k) = (\theta_k-\theta_d)^2 \;+\; 0.01\,(\omega_k-\omega_d)^2,
%\label{eq:pendulum_cost} \\
%c(\bv{u}_k) = 0.1u^2,
%\end{gather}
%with $\bv{x}_k = [\theta_k,\ %\omega_k]^\top$.

%The discrete-time stochastic dynamics used to generate data are
%\begin{equation}
%\begin{aligned}
%\theta_k &= \theta_{k-1} + \omega_{k-1}\,dt + W_\theta,\\
%\omega_k &= \omega_{k-1} + \Big(\tfrac{g}{l}\sin(\theta_{k-1}) + \tfrac{u_k}{m\,l^2}\Big)\,dt + W_\omega,
%\end{aligned}
%\label{eq:pendulum_dyn}
%\end{equation}
%where $dt = ??? \,\mathrm{s}$, $W_\theta \sim \mathcal{N}(0,\,0.05)$ and $W_\omega \sim \mathcal{N}(0,\,0.1)$.
%The control input satisfies $u_k \in U \coloneqq [-2.5,2.5]$.

%In the running example, we let $X_k := [\theta_k,\ \omega_k]^\top$ and consider two pendulums (``target'' and ``reference'') with different physical parameters. For the \textit{target} pendulum (the one to be controlled), the parameters are $m=1\,\mathrm{kg}$ and $l=0.6\,\mathrm{m}$. For the \textit{reference} pendulum, the parameters are $m=0.5\,\mathrm{kg}$ and $l=0.5\,\mathrm{m}$.

%To construct probabilistic models used in the paper, the state domain was taken as $X := [-\pi,\pi]\times[-5,5]$.}

\addtolength{\textheight}{-4cm} 
\section{Conclusion}

We considered the problem of computing optimal policies under a free energy minimization framework that accounts for ambiguity in the generative model. This leads to a distributionally robust minimax formulation of the control problem. We showed that the optimal policy can be computed via a bi-level optimization approach, where the free energy is first maximized within the ambiguity set and then minimized in the policy space.
We then showed that the inner maximization can be recast as a DC program, and we adapted a Frank–Wolfe algorithm, providing  convergence guarantees, to tackle the  optimization. The effectiveness of the results was illustrated on a pendulum stabilization task. Based on these early results, future work includes: (i)  embedding learning in policy computation; (ii) extending our results to include ambiguity on $\refpolicy{k}{k-1}$.

\noindent{\bf Acknowledgments.} We acknowledge the use of ChatGPT for assistance in improving wording and grammar of this document.


%%%%%%%%%%%%%%%%%%%%%%%%%%%%%%%%%%%%%%%%%%%%%%%%%%%%%%%%%%%%%%%%%%%%%%%%%%%%%%%%

%%%%%%%%%%%%%%%%%%%%%%%%%%%%%%%%%%%%%%%%%%%%%%%%%%%%%%%%%%%%%%%%%%%%%%%%%%%%%%%%

%%%%%%%%%%%%%%%%%%%%%%%%%%%%%%%%%%%%%%%%%%%%%%%%%%%%%%%%%%%%%%%%%%%%%%%%%%%%%%%%
%\section*{APPENDIX}

%\section*{Acknowledgment}
%Caio César Graciani Rodrigues gratefully acknowledges funding from the Scuola Superiore Meridionale, MERC department.

\bibliographystyle{abbrv}        % Include this if you use bibtex 
\bibliography{russobib-2}           % and a bib file to produce the 

\end{document}